\documentclass{amsart}

\usepackage{amsthm, amssymb, amsmath}
\usepackage[english]{babel}
\usepackage{graphicx}
\usepackage{subfigure}

\begin{document}

\title[A natural extension for the greedy $\beta$-transformation]{A natural extension for the greedy $\beta$-transformation with three deleted digits}
\author{Karma Dajani}
\address{Department of Mathematics\\
Utrecht University\\
Postbus 80.000\\
3508 TA Utrecht\\
the Netherlands} \email{k.dajani1@uu.nl}
\author{Charlene Kalle}
\address{Department of Mathematics\\
Utrecht University\\
Postbus 80.000\\
3508 TA Utrecht\\
the Netherlands} \email{c.c.c.j.kalle@uu.nl}

\maketitle

\begin{abstract}
We give an explicit expression for the invariant measure, absolutely continuous with respect to the Lebesgue measure, of the greedy $\beta$-transformation with three deleted digits. We define a version of the natural extension of the transformation to obtain this expression. We get that the transformation is exact and weakly Bernoulli.
\end{abstract}

\newtheorem{prop}{Proposition}[section]
\newtheorem{theorem}{Theorem}[section]
\newtheorem{lemma}{Lemma}[section]
\newtheorem{cor}{Corollary}[section]
\newtheorem{remark}{Remark}[section]
\theoremstyle{definition}
\newtheorem{defn}{Definition}[section]
\newtheorem{ex}{Example}[section]
\newcommand{\lex}{<_{\text{lex}}}
\newcommand{\T}{\mathcal T}

\section{Introduction}

A $\beta$-expansion with deleted digits of a number $x$ is an expression of the form
\begin{equation}\label{q:exp}
x = \sum_{n=1}^{\infty} \frac{b_n}{\beta^n},
\end{equation}
where $\beta >1$ is a real number and the $b_n$ are all elements of a set of real numbers $A = \{ a_0, a_1, \ldots, a_m \}$. In \cite{Ped1}, Pedicini showed that if $\beta >1$ and if the set of real numbers $A=\{ a_0, a_1, \ldots, a_m \}$ satisfies the following two conditions:
\begin{itemize}
\item[(i)] $a_0 < a_1 < \ldots < a_m$,
\item[(ii)] $\displaystyle \max_{1\leq j \leq m} (a_j - a_{j-1}) \leq \frac{a_m - a_0}{\beta -1}$,
\end{itemize}
then every $x \in \displaystyle \left[ \frac{a_0}{\beta -1}, \frac{a_m}{\beta -1} \right]$ has an expansion of the form (\ref{q:exp}), with $b_n \in A$ for all $n \ge 1$. He also gave an iterative algorithm that generates greedy expansions of the form (\ref{q:exp}), in the sense that at each step of the algorithm, $b_n$ is chosen to be the largest element of $A$ such that
\begin{equation}\label{q:greedy}
\sum_{i=1}^n \frac{b_i}{\beta^i} + \sum_{i=n+1}^{\infty} \frac{a_1}{\beta^i} \le x.
\end{equation}
In \cite{DK2} a dynamical system is given that generates all possible expansions of the form (\ref{q:exp}) for points $x \in \displaystyle \left[ \frac{a_0}{\beta -1}, \frac{a_m}{\beta -1} \right]$, if $A$ satisfies (i) and (ii). We call a set $A$, satisfying (i) and (ii) an {\it allowable digit set}. If $A$ is not allowable, then not every point has an expansion. The size of the set of real numbers that can be represented by expression (\ref{q:exp}), for different choices of $\beta$ and for $A = \{ 0,1,3 \}$, was already studied by Keane, Smorodinsky and Solomyak in \cite{Kea1}. Pollicott and Simon studied the Hausdorff dimension of this set of reals in the case $A$ is a subset of the non-negative integers in \cite{Pol1}.\\
\indent The $\beta$-expansions with deleted digits are a generalization of the $\beta$-expansions with what we could call a complete digit set, i.e. with $A = \{ 0,1, \ldots , \lfloor \beta \rfloor \}$, where $\lfloor x \rfloor$ denotes the largest integer less than or equal to $x$. A lot of research has been done on this topic and we mention some results here. For a given $\beta$, all $x \in \displaystyle \left[ 0, \frac{\lfloor \beta \rfloor}{\beta -1} \right]$ have an expansion with a complete digit set and almost all $x$ in this interval have continuum expansions of this form. (See \cite{Erd1} by Erd\"os, Jo\'o and Komornik and \cite{Sid1} by Sidorov for more information.) One way to generate $\beta$-expansions with a complete digit set is by iterating the transformation
$$ T_c : \left[ 0, \frac{\lfloor \beta \rfloor}{\beta -1} \right] \to  \left[ 0, \frac{\lfloor \beta \rfloor}{\beta -1} \right]:  x \to \left\{
\begin{array}{ll}
\beta x \; (\text{mod }1), & \text{if } 0 \leq x < \displaystyle \frac{\lfloor \beta \rfloor}{\beta},\\
\\
\beta x - \lfloor \beta \rfloor, & \text{if } \displaystyle  \frac{\lfloor \beta \rfloor}{\beta} \leq x \leq  \frac{\lfloor \beta \rfloor}{\beta-1}.
\end{array}
\right.
$$
If we set
\begin{equation}\label{q:digitset}
b^c_1 =b^c_1(x) = \left\{
\begin{array}{ll}
i, & \mbox{if } x \in \left[ \displaystyle\frac{i-1}{\beta}, \frac{i}{\beta} \right), \mbox{ for } i=1, \ldots, \lfloor \beta \rfloor,\\
\\
\lfloor \beta \rfloor, & \mbox{if } x \in \left[\displaystyle\frac{\lfloor \beta \rfloor}{\beta}, \frac{\lfloor \beta \rfloor}{\beta-1} \right],
\end{array}
\right.
\end{equation}
and for $n \ge 1$, $b^c_n = b^c_n(x) = b^c_1(T_c^n x)$, then $T_c x=\beta x-b^c_1$, and for
any $n\ge 1$,
$$x=\sum_{i=1}^n\frac{b^c_i}{\beta^i}+\frac{T_c^n x}{\beta^n}.$$
Letting $n\to \infty$, it is easily seen that
$x=\displaystyle\sum_{n=1}^{\infty}\frac{b^c_n}{\beta^n}.$
We call the $\beta$-expansion generated by the transformation $T_c$ the greedy $\beta$-expansion with a complete digit set of $x$, since it gives us at each step in the expansion the largest digit possible. More precisely, for each $n \ge 1$, if $b^c_1, \ldots, b^c_{n-1}$ are already known, then $b^c_n$ is the largest element of the complete digit set $A = \{ 0, 1, \ldots \lfloor \beta \rfloor \}$, such that
$$ \sum_{i=1}^{n} \frac{b^c_i}{\beta ^i} \le x. $$
Throughout the rest of the paper, $\lambda$ will denote the 1-dimensional Lebesgue measure. The transformation $T_c$ has a unique invariant measure, absolutely continuous with respect to $\lambda$. R\'enyi proved the existence of this measure (\cite{Ren1}) and Gel'fond and Parry independently gave an explicit formula for the density function of this measure in \cite{Gel1} and \cite{Par1} respectively. The invariant measure has the unit interval $[0,1)$ as its support and the density function $h_c$ is given as follows.
\begin{equation}\label{q:parrymeasure}
h_c: [0,1) \to [0,1) : x \mapsto \frac{1}{F(\beta)} \sum_{n=0}^{\infty} \frac{1}{\beta^n} 1_{[0, T_c^n 1)}(x),
\end{equation}
where $F(\beta) = \displaystyle \int_0^1 \sum_{x < T_c^n 1} \frac{1}{\beta^n} d\lambda$ is a normalizing constant. From now on we will refer to an invariant measure, absolutely continuous with respect to $\lambda$ as an {\it acim}.\\
\indent For the $\beta$-expansions with deleted digits we also have a transformation that generates the greedy expansions by iteration. This transformation is called the greedy $\beta$-transformation with deleted digits. If $\beta >1$ is a real number and $A= \{ 0, a_1, \ldots, a_m \}$ is an allowable digit set of which the first digit equals zero, i.e. $a_0=0$, then this transformation $T: \displaystyle \left[0, \frac{a_m}{\beta -1} \right] \to \left[0, \frac{a_m}{\beta -1} \right]$ is given by
$$ Tx = \left\{
\begin{array}{ll}
\beta x - a_j, & \mbox{if } x \in \left[ \displaystyle \frac{a_j}{\beta}, \frac{a_{j+1}}{\beta} \right), \mbox{ for } j=0, \ldots, m-1,\\
\\
\beta x - a_m, & \mbox{if } x \in \left[ \displaystyle\frac{a_m}{\beta}, \frac{a_m}{\beta -1} \right]. 
\end{array}
\right.
$$
In \cite{DK2} it is shown that a greedy $\beta$-transformation with deleted digits for a given $\beta >1$ and an allowable digit set $A=\{ a_0, \ldots, a_m \}$ is isomorphic to a greedy $\beta$-transformation with deleted digits for the same $\beta >1$, but with allowable digit set $\tilde{A} = \{ 0, a_1-a_0, \ldots, a_m-a_0 \}$, i.e. a digit set of which the first digit equals zero. Therefore we can assume without loss of generality that $a_0=0$ for all allowable digit sets. The sequence of digits $\{ b_n \}_{n \ge 1}$ can be defined in a way similar to (\ref{q:digitset}) as follows. Put $b_1(x) = a_j$ if $x \in \left[\displaystyle \frac{a_j}{\beta}, \frac{a_{j+1}}{\beta}\right)$ and $b_1(x) = a_m$ if $x \in \left[ \displaystyle\frac{a_m}{\beta}, \frac{a_m}{\beta -1} \right]$ and for $n \ge 1$, put $b_n (x) = b_1(T^{n-1}x)$. The definition of the transformation $T$ is based on the greedy algorithm with deleted digits that was defined by Pedicini in \cite{Ped1} and the transformation therefore generates greedy $\beta$-expansions with deleted digits satisfying (\ref{q:greedy}). (For more information, see \cite{DK2}.) If $x = \sum_{n=1}^{\infty} \frac{b_n}{\beta ^n}$ is the greedy $\beta$-expansion with deleted digits for $x$, we also write
$$x =_{\beta} b_1 b_2 b_3 \ldots,$$
which is understood to mean the same as (\ref{q:exp}).\\
\indent From \cite{DK1} we know that the transformation $T$ has an acim that is unique and ergodic. The support of this measure is given by the interval $[0, a_{i_0}-a_{i_0-1})$, where
\begin{equation}\label{q:i0}
i_0 = \min \{i \in \{ 1,\ldots ,m\} : T[0,a_i-a_{i-1}) \subseteq [0, a_i-a_{i-1}) \; \lambda \text{ a.e. } \}.
\end{equation}
An explicit expression for the density function of this measure however, is given only for special cases. One of these special cases is, when 
\begin{equation}\label{q:critwilkinson}
m < \beta \le m+1.
\end{equation}
Before we can give this density, however, we need some more notation. In \cite{DK2} it is proven that the minimal amount of digits in an allowable digit set is $\lceil \beta \rceil$. In other words, the amount of digits in $A$ is at least equal to the smallest integer larger than or equal to $\beta$.\\
\indent Let $N$ be the largest element of the set $\{ 1, \ldots, m \}$ such that $\displaystyle \frac{a_N}{\beta}< a_{i_0}-a_{i_0-1}$. Define a partition $\Delta = \{ \Delta (a_i) : 0 \le i \le N \}$ of the support of the acim of $T$, where for $i=0, \ldots, N-1$, we have
$$ \Delta (a_i) = \left[\frac{a_i}{\beta}, \frac{a_{i+1}}{\beta} \right),$$
and
$$ \Delta(a_N) = \left[\frac{a_N}{\beta},a_{i_0}-a_{i_0-1} \right).$$
Note that $T \Delta(a_N) = [0, T(a_{i_0}-a_{i_0-1}))$ and for $i \in \{0, \ldots, N-1\}$, $T \Delta(a_i) = [0,a_{i+i}-a_i)$. Using $\Delta$ and $T$, we can make a sequence of partitions $\{ \Delta^{(n)} \}$: for $n \geq 0$, 
\begin{equation}\label{q:deltan}
\displaystyle \Delta^{(n)} = \bigvee_{k=0}^{n-1} T^{-k} \Delta.
\end{equation}
The elements of $\Delta^{(n)}$ are intervals and are called the {\it fundamental intervals of rank} $n$. If
$$ \Delta(b_0) \cap T^{-1} \Delta(b_1) \cap \ldots \cap T^{-(n-1)} \Delta(b_{n-1})$$ 
is an element of $\Delta^{(n)}$, denote it by $\Delta (b_0 b_1 \ldots b_{n-1})$. We call a fundamental interval $\Delta (b_0 b_1 \ldots b_{n-1})$ {\it full} of rank $n$ if 
$$\lambda (T^n \Delta (b_0 b_1 \ldots b_n)) = a_{i_0}-a_{i_0-1}$$
and {\it non-full} otherwise. This means that for a full fundamental interval, $\Delta (b_0 \ldots b_{n-1})$ we have
\begin{equation}\label{q:genfull}
\lambda (\Delta (b_0 b_1 \ldots b_{n-1})) = \frac{a_{i_0}-a_{i_0-1}}{\beta^n}
\end{equation}
and if $\Delta (b_0 \ldots b_{n-1})$ is non-full, then
\begin{equation}\label{q:gennonfull}
\lambda (\Delta (b_0 b_1 \ldots b_{n-1})) < \frac{a_{i_0}-a_{i_0-1}}{\beta^n}.
\end{equation}
Let $B_n$ be the collection of all non-full fundamental intervals of rank $n$, that are not subsets of any full fundamental interval of lower rank. For $x \in [0, a_{i_0}-a_{i_0-1})$, define $\phi_0 (x) =1$ and for $n \geq 1$, let 
$$ \phi_n(x) = \sum_{\Delta (b_0 b_1 \ldots b_{n-1}) \in B_n} \frac{1}{\beta^n} 1_{T^n \Delta (b_0 b_1 \ldots b_{n-1})}(x).$$
Put $\phi = \displaystyle \sum_{n=0}^{\infty} \phi_n$. Then the function
\begin{equation}\label{q:wilkinson}
h:[0, a_{i_0}-a_{i_0-1}) \rightarrow [0, a_{i_0}-a_{i_0-1}): x \mapsto \displaystyle \frac{\phi(x)}{ \int \phi(x) d\lambda(x)}
\end{equation}
is the density function of the acim of $T$. (This density is a special case of the density found by Wilkinson in \cite{Wil1}.) Notice that for the classical greedy $\beta$-transformation, $B_n$ contains at most one element $\Delta(b_0 \ldots b_{n-1})$ and for this element we have
$$ T_c^n \Delta(b_0 \ldots b_{n-1}) = [0, T^n_c 1).$$
So the density function from (\ref{q:parrymeasure}) is a special case of the density function from (\ref{q:wilkinson}).\\
\indent In this article we will give an expression for the density function of the acim of the greedy $\beta$-transformation with three deleted digits. This means that we will only be looking at $1 < \beta < 3$ and allowable digit sets $A$ of the form $A = \{ 0, a_1, a_2 \}$. Since $2 < \beta <3$ fits the framework above, we already know that the density is given by (\ref{q:wilkinson}). However, (\ref{q:wilkinson}) will also turn out to be valid for $1 < \beta \le 2$. In the first section we will introduce some notation and discuss when we can directly relate the density function of $T$ to the density from equation (\ref{q:wilkinson}). For the case in which this is not immediately clear, we will derive a formula for the density by defining a version of the natural extension of the dynamical system
$$([0, a_{i_0}-a_{i_0-1}), \mathcal B ([0, a_{i_0}-a_{i_0-1})), \mu, T).$$
Here $\mathcal B ([0, a_{i_0}-a_{i_0-1}))$ is the Borel $\sigma$-algebra on the support of the acim of $T$ and $\mu$ is the measure obtained by ``pulling back'' the appropriate measure, defined on the natural extension. We will define this version of the natural extension in the second section. The density of the measure $\mu$ will again be the density from (\ref{q:wilkinson}). For the greedy $\beta$-transformation with a complete digit set, a version of the natural extension is given in \cite{DKS} and also in \cite{Bro1} by Brown and Yin. Our definition is based on the version given in the first paper. Using the measure $\mu$, we obtain that the transformation $T$ is exact and weakly Bernoulli. To illustrate the construction of our version of the natural extension, an example of a specific greedy $\beta$-transformation with three deleted digits can be found in the last section. Here $\beta$ is the golden mean and $A = \{ 0,3,4\}$.

\section{A closer look at the greedy $\beta$-transformation with three deleted digits}
\begin{figure}[h]
\centering
\subfigure[$\beta \; = \; \sqrt 3 \;$ and $ A=\{ 0,1,3 \}$]{\includegraphics[height=2.8cm]{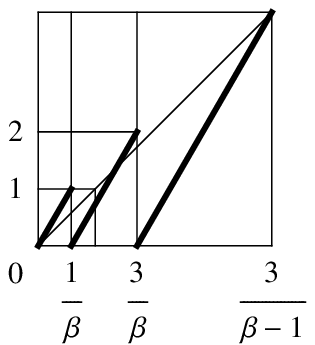}}
\quad
\subfigure[$\beta = 1 + \sqrt 2$ and $A=\{ 0,1,3\}$]{\includegraphics[height=2.8cm]{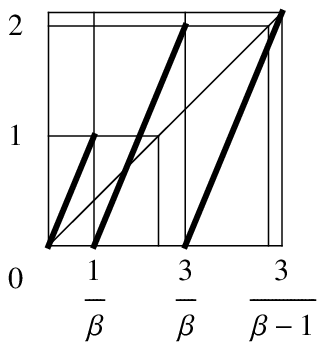}}
\quad
\subfigure[$\beta \; = \;  \sqrt 7 \;$ and $A=\{ 0,3,7\}$]{\includegraphics[height=2.8cm]{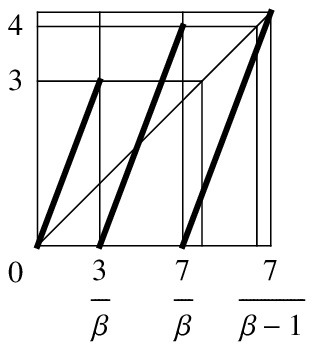}}
\quad
\subfigure[$\beta = \frac{1+\sqrt 5}{2} $ and $A=\{ 0,3,4 \}$]{\includegraphics[height=2.8cm]{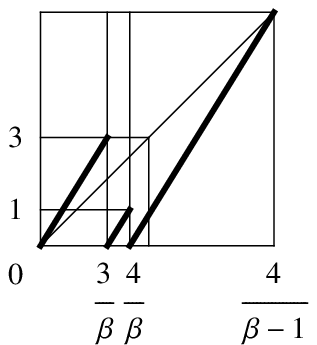}}
\caption{Examples of the four possibilities for the support of the acim of the greedy $\beta$-transformation with three deleted digits.}
\label{f:attractorcases}
\end{figure}

The aim of this paper is to give an explicit expression for the acim of the greedy $\beta$-transformation with three deleted digits. For $2 < \beta <3$, we already know that the density function of this measure is given by (\ref{q:wilkinson}). We will first give some conditions under which we can directly see that the density from (\ref{q:wilkinson}) is the density of the acim.\\
\indent Let $T$ be the greedy $\beta$-transformation with deleted digits for which $1 < \beta < 3$ and for which $A = \{ 0, a_1, a_2 \}$ is an allowable digit set. So, we set $a_0=0$. The transformation then becomes:
$$ T x = \left\{
\begin{array}{ll}
\beta x, & \text{if } x \in \left[ 0,\displaystyle  \frac{a_1}{\beta} \right),\\
\beta x - a_1, & \text{if } x \in \left[ \displaystyle \frac{a_1}{\beta}, \frac{a_2}{\beta} \right),\\
\beta x -a_2, & \text{if } x \in \left[ \displaystyle \frac{a_2}{\beta}, \frac{a_2}{\beta -1} \right].
\end{array}
\right.
$$
The sequence of digits $\{b_n \}_{n \ge 1}$ is defined for any $x \in \left[ \displaystyle 0, \frac{a_2}{\beta-1} \right]$ as in the introduction. Before we can give the density function of the acim, we need to find its support. From (\ref{q:i0}) we know that this support is either the interval $[0, a_1)$ or the interval $[0, a_2-a_1)$. The following situations can occur.\\
\indent  Suppose first that $a_1 < \frac{a_2}{\beta}$. Either $Ta_1 \le a_1$ or $T a_1 > a_1$.\\
\indent \indent If $Ta_1 \le a_1$, then the support is $[0, a_1)$ and since $T a_1 = \beta a_1 -a_1$, we get that $\beta \le 2$. The transformation $T$ on the interval $[0, a_1)$ is then isomorphic to the classical greedy $\beta$-transformation, $T_c$, for the same $\beta$. Figure \ref{f:attractorcases}(a) is an example of this.\\
\indent \indent If $T a_1 > a_1$, then the support of the invariant measure is the other interval, $[0, a_2-a_1)$, and we can deduce that $\beta >2$. So, criterion (\ref{q:critwilkinson}) applies and the density for the invariant measure is given by (\ref{q:wilkinson}). In Figure \ref{f:attractorcases}(b) we see an example.\\
\indent Suppose that $a_1 \ge \frac{a_2}{\beta}$. Either $a_2-a_1 > a_1$ or $a_2 -a_1 \le a_1$.\\
\indent \indent If $a_2 - a_1 > a_1$, then the support of the acim is $[0, a_2-a_1)$ and we have that $2a_1 < a_2 \le \beta a_1$. So again $\beta >2$ and the density from equation (\ref{q:wilkinson}) is the density of the invariant measure we are looking for. See Figure \ref{f:attractorcases}(c) for an example.\\
\indent \indent If $a_2-a_1 = a_1$, then the support is $[0, a_1)=[0, a_2-a_1)$ and $\beta \ge 2$. The transformation in this case is isomorphic to the classical greedy transformation for the same $\beta$. Similarly, if $a_2 -a_1 < a_1$ and $a_1 = \frac{a_2}{\beta}$, then $\beta \le 2$ and the support of the acim is $[0, a_1)$. The transformation $T$ is again isomorphic to the classical greedy $\beta$-transformation for the same $\beta$.\\
\indent \indent Lastly, suppose that $a_1 = \frac{a_2}{\beta -1}$. If $\beta \le 2$, then this would imply
$$a_2 = (\beta -1) a_1 \le a_1,$$
which gives a contradiction. So this means that $\beta >2$ and thus that the density from (\ref{q:wilkinson}) is the density for the acim.\\
\indent The only situation we did not yet consider, is when $ \frac{a_2}{\beta} < a_1 < \frac{a_2}{\beta -1}$ and $a_2-a_1 < a_1$. To this, the rest of the paper is dedicated. Figure \ref{f:attractorcases}(d) gives an example of a transformation that satisfies these conditions.

From now on, suppose that $T$ is a greedy $\beta$-transformation with an allowable digit set $A=\{ 0, a_1, a_2\}$, that satisfies the following equations: $a_1 > a_2-a_1$ and $a_1 > \displaystyle \frac{a_2}{\beta}$. Of course, since $A$ is allowable and $(\beta -1)a_1 \neq a_2$, we also have $a_1 < \displaystyle \frac{a_2}{\beta -1}$. This is all captured in the following condition.
\begin{equation}\label{q:conditions}
a_1 \cdot \max \{ \beta -1, 1 \} < a_2 < a_1 \cdot \min \{2, \beta \}.
\end{equation}
Notice that we do not assume that $\beta \le 2$, although we already know that the density from (\ref{q:wilkinson}) is the density of the acim in this case. The reason is that the construction of the natural extension that will be given in what follows, is also valid for $2 < \beta < 3$.

\begin{remark}\label{r:1beta2}
{\rm Observe that, if $1 < \beta \le 2$, then $\displaystyle a_1 > \frac{a_2}{\beta}$ implies condition (\ref{q:conditions}). So, if $1 < \beta \le 2$ and a digit set $A= \{ 0,a_1, a_2 \}$ satisfies $\beta a_1 > a_2$, then $A$ is an allowable digit set.
}
\end{remark}

The support of the acim of the transformation $T$ is the interval $[0, a_1)$. The partition $\Delta = \{ \Delta (0), \Delta (a_1), \Delta (a_2) \}$ of this interval is given in the following way:
$$\Delta (0) = \left[ 0, \frac{a_1}{\beta} \right), \quad \Delta (a_1) = \left[ \frac{a_1}{\beta}, \frac{a_2}{\beta} \right), \quad \Delta (a_2) = \left[ \frac{a_2}{\beta}, a_1 \right).$$
As explained in the introduction, we can construct the sequence of partitions $\{ \Delta^{(n)} \}_{n \geq 1}$, with $\Delta^{(n)}$ as defined in (\ref{q:deltan}). The elements of $\Delta^{(n)}$ are denoted by $\Delta(b_0 \ldots b_{n-1})$ and are either full or non-full fundamental intervals of rank $n$. We have the following obvious lemma.

\begin{lemma}\label{l:fullfi}
Let $\Delta(d_0 \ldots d_{p-1})$ and $\Delta(e_0 \ldots e_{q-1})$ be two full fundamental intervals of rank $p$ and $q$ respectively. Then $\Delta (d_0 \ldots d_{p-1}e_0 \ldots e_{q-1})$ is a full fundamental interval of rank $p+q$.
\end{lemma}

Recall that for $n \ge 1$, $B_n$ is the collection of all non-full fundamental intervals of rank $n$ that are not contained in any full fundamental interval of lower rank. Let $\kappa (n)$ be the number of elements in $B_n$. So $\kappa (1) =2$, since this is the number of non-full fundamental intervals of rank 1 and for all $n \ge 1$, $\kappa (n) \le 2^n$. The version of the natural extension that we will define in the next section, uses all the elements of $B_n$ for all $n \ge 1$. To make sure that the total measure of the underlying space of this version is finite, we need upper bounds for the values of $\kappa (n)$. To obtain these, we will first describe the structure of the elements of $B_n$.\\
\indent Notice that by (\ref{q:gennonfull}) we have that for all elements $\Delta(b_0 \ldots b_{n-1}) \in B_n$,
\begin{equation}\label{q:nonfull}
\lambda (\Delta(b_0 \ldots b_{n-1})) < \frac{a_1}{\beta^n}.
\end{equation}
For an $x \in [0, a_1)$, the set $\{ T^n x : n \ge 0 \}$ is called the {\it orbit} of $x$ under the transformation $T$. Let
\begin{eqnarray*}
a_2 - a_1 &=& \sum_{n=1}^{\infty} \frac{d^{(1)}_n}{\beta^n} =_{\beta} d^{(1)}_1d^{(1)}_2 d^{(1)}_3 \ldots ,\\
\beta a_1 - a_2 &=& \sum_{n=1}^{\infty} \frac{d^{(2)}_n}{\beta^{n+1}} =_{\beta} d^{(2)}_1 d^{(2)}_2 d^{(2)}_3 \ldots ,
\end{eqnarray*}
be the greedy $\beta$-expansions with deleted digit set $A$ of the points $a_2 - a_1$ and $\beta a_1 - a_2$, i.e. the expansions generated by iterations of $T$. The number $\beta a_1 - a_2$ would be the image of $a_1$ under $T$ if $T$ were not restricted to the interval $[0,a_1)$. The values of the numbers $\kappa (n)$ depend on the orbits of the points $\beta a_1-a_2$ and $a_2 -a_1$ under $T$ and whether or not $T^i (a_2-a_1)$ and $T^i (\beta a_1-a_2)$ are elements of $\Delta(a_2)$ for $0 \le i < n$. To see this, notice that for any $\Delta(b_0 \ldots b_{n-1}) \in B_n$ one has $b_0 \in \{ a_1, a_2 \}$ and the set $T^n \Delta(b_0 \ldots b_{n-1})$ has the form
$$[0, T^i (a_2-a_1)) \text{ or } [0, T^i (\beta a_1 - a_2))$$
for some $0 \le i <n$. Suppose $T^n \Delta(b_0 \ldots b_{n-1}) = [0, T^i (a_2-a_1))$.\\
\indent If $\lambda (T^i (a_2 -a_1) \cap \Delta(a_2))=0$, then $\Delta(b_0 \ldots b_{n-1})$ contains exactly one element of $B_{n+1}$, namely $ \Delta(b_0 \ldots b_{n-1} 0)$ in case $\lambda (T^i (a_2 -a_1) \cap \Delta(a_1)) =0$ or $ \Delta(b_0 \ldots b_{n-1} a_1)$ in case $\lambda (T^i (a_2 -a_1) \cap \Delta(a_1))>0$. Furthermore, in the first case,
$$T^{n+1} \Delta(b_0 \ldots b_{n-1} 0) = [0, T^{i+1} (a_2-a_1))$$
and also in the second case,
$$T^{n+1} \Delta(b_0 \ldots b_{n-1} a_1) =  [0, T^{i+1} (a_2-a_1)).$$
\indent On the other hand, if $\lambda (T^i (a_2 - a_1) \cap \Delta(a_2))>0$, then $\Delta(b_0 \ldots b_{n-1})$ contains exactly two elements of $B_{n+1}$, namely the sets $\Delta(b_0 \ldots b_{n-1} a_1)$ and $\Delta(b_0 \ldots b_{n-1} a_2)$. Now,
$$T^{n+1} \Delta(b_0 \ldots b_{n-1} a_1) = [0, a_2-a_1)$$
and
$$T^{n+1} \Delta(b_0 \ldots b_{n-1} a_2) = [0, T^{i+1} (a_2-a_1)).$$
Similar arguments hold in case $T^n \Delta(b_0 \ldots b_{n-1}) = [0, T^i (\beta a_1-a_2))$, except that $T^i (a_2-a_1)$ is replaced by $T^i (\beta a_1-a_2)$.\\
\indent For $n \ge 1$, let $\bar{\kappa} (n)$ be the number of elements from $B_n$ that contain two elements from $B_{n+1}$. Then clearly for all $n \ge 1$,
\begin{equation}\label{q:barkappa}
\kappa (n+1) = \kappa(n) + \bar{\kappa}(n).
\end{equation}
From the above we see that in order to get an upper bound on $\kappa (n)$, we need to study the behavior of the orbits of $a_2-a_1$ and $\beta a_1-a_2$. The following three lemmas say something about the first few elements of the orbits of these points. These lemmas are needed to guarantee the total measure of the underlying space of the natural extension will be finite.

\begin{lemma}\label{l:Gbeta2}
If $1 < \beta \le 2$ and $\displaystyle a_1 > \frac{a_2}{\beta}$, then $a_2-a_1 \not \in \Delta (a_2)$.
\end{lemma}

\begin{proof}
Since $\beta \le 2$, we have $ 1-\frac{1}{\beta} \le \frac{1}{\beta}$. Thus $ a_2 \, \left( 1-\frac{1}{\beta} \right) \le \frac{a_2}{\beta} < a_1 $
and hence $a_2 - a_1 < \frac{a_2}{\beta}$. This proves the lemma.
\end{proof}

Observe that $B_n$ only contains fundamental intervals of which the first digit is either $a_1$ or $a_2$. Let $\kappa_1(n)$ denote the number of elements $\Delta (b_0 \ldots b_{n-1})$ in $B_n$ such that $b_0 = a_1$ and $\kappa_2(n)$ the number of elements in $B_n$ that have $a_2$ as their first digit. Then of course for all $n \ge 1$,
$$\kappa(n) = \kappa_1(n) + \kappa_2(n).$$
Let $\{F(n) \}_{n \ge 0}$ denote the Fibonacci sequence, i.e. let $F(0)=0$, $F(1)=1$ and for $n \ge 2$, let $F(n) = F(n-1)+F(n-2)$. Lemma \ref{l:Gbeta2} implies that the number of elements of $B_{n+1}$ would be maximal if the only elements of $B_n$ that do not contain two elements from $B_{n+1}$ are the elements $\Delta(b_0 \ldots b_{n-1})$ for which $T^n \Delta(b_0 \ldots b_{n-1}) =a_2-a_1$. In this maximal situation we would have $\kappa_1(1) = \kappa_1 (2) =1$ and for $n \ge 3$, 
$$\kappa_1(n) = \kappa_1(n-1) + \kappa_1(n-2).$$
For $\kappa_2$ we would have that $\kappa_2(n) = \kappa_1(n+1)$. This means that under the conditions from Lemma \ref{l:Gbeta2}, we have for all $n \ge 1$ that $\kappa_1 (n) \le F(n)$ and
\begin{equation}\label{q:Gbeta2}
\kappa (n)= \kappa_1(n)+ \kappa_2(n) \le F(n)+F(n+1)=F(n+2).
\end{equation}

\vskip .3cm

\noindent Let $G = \frac{1+ \sqrt 5}{2}$ be the golden mean, i.e. the positive solution of the equation $x^2 -x-1=0$.
\begin{lemma}\label{l:lessthanG}
Let $1 < \beta \le G$ and $a_1 > \displaystyle \frac{a_2}{\beta}$. Then $a_2-a_1, \, \beta a_1 - a_2 \in \Delta (0)$.
\end{lemma}

\begin{proof}
Since $\beta \le G$, we have $1+ \frac{1}{\beta} \ge \beta$, so by equation (\ref{q:conditions}),
$ \left( 1+ \frac{1}{\beta} \right) \, a_1 \ge \beta a_1 > a_2$.
Thus $a_2 - a_1 <  \frac{a_1}{\beta}$ and hence $a_2 - a_1 \in \Delta (0)$. On the other hand,
since $\beta -  \frac{1}{\beta} \le 1$, we have $ \left( \beta - \frac{1}{\beta} \right) \, a_1 \le a_1 < a_2$.
Thus $\beta a_1 - a_2 <  \frac{a_1}{\beta}$ and we have that $\beta a_1 -a_2 \in \Delta (0)$.
\end{proof}

\begin{remark}\label{r:wortel2}
{\rm This lemma implies that for $1 < \beta \le G$ and for digit sets satisfying condition (\ref{q:conditions}),  we have $\kappa (2)=2$. The largest amount of elements for $B_n$ would be obtained if $\lambda (T^i (a_2-a_1) \cap \Delta(a_2)) >0$ and $\lambda(T^i(\beta a_1-a_2) \cap \Delta(a_2))>0$ for all odd values of $i$ and thus $T^i (a_2-a_1), T^i(\beta a_1-a_2) \in \Delta(0)$ for all even values of $i$. In this case,
$$
\bar{\kappa}(n) = \left\{
\begin{array}{ll}
\kappa(n), & \text{if $n$ is odd},\\
0, & \text{if $n$ is even}.
\end{array}
\right.
$$
This would imply that for $n \ge 1$, $\kappa (2n-1) = \kappa(2n) = 2^n$. So, in general we have that for all $n \ge 1$, $\kappa (n) \le 2^{\lfloor n/ 2 \rfloor +1}$.
}
\end{remark}

\begin{lemma}\label{l:2to1overm}
Let $m \ge 2$ and $a_1 > \displaystyle \frac{a_2}{\beta}$. If $1 < \beta \le 2^{\frac{1}{m}}$, then 
$$T^i (a_2 - a_1), \, T^i (\beta a_1-a_2) \in \Delta (0)$$
for all $i \in \{ 0, 1, \ldots, m-1 \}$.
\end{lemma}

\begin{proof}
The proof is by induction on $m$. Note that from Lemma \ref{l:lessthanG} we know that $a_2-a_1, \, \beta a_1-a_2 \in \Delta (0)$. Assume first that $m=2$ and thus $\beta \le \sqrt 2 = 2^{\frac{1}{2}}$. Then
\begin{eqnarray*}
T(a_2-a_1) = \beta (a_2-a_1) \in \Delta (0) & \Leftrightarrow & \beta (a_2-a_1) < \frac{a_1}{\beta}\\
& \Leftrightarrow & \frac{a_2}{\beta} < a_1 \, \left( \frac{1}{\beta} + \frac{1}{\beta^3} \right).
\end{eqnarray*}
If $\beta \le \sqrt 2$, then $\displaystyle \frac{1}{\beta} + \frac{1}{\beta^3}>1$ and $\displaystyle \frac{a_2}{\beta} < a_1 < \left( \frac{1}{\beta}+ \frac{1}{\beta^3} \right) \, a_1$. So $T(a_2 - a_1) \in \Delta (0)$. On the other hand, since $\beta a_1-a_2 \in \Delta (0)$, we have
\begin{eqnarray*}
T(\beta a_1-a_2) = \beta (\beta a_1 - a_2)\in \Delta (0) &\Leftrightarrow & \beta^2 a_1 - \beta a_2 < \frac{a_1}{\beta}\\
&\Leftrightarrow & \left( \beta - \frac{1}{\beta^2} \right) \, a_1 < a_2.
\end{eqnarray*}
If $\beta \le \sqrt 2$, then $\displaystyle \beta -\frac{1}{\beta^2} < 1$. So $\displaystyle \left( \beta - \frac{1}{\beta^2} \right) \, a_1 < a_1 < a_2$ and thus $T(\beta a_1-a_2) \in \Delta (0)$.\\
Now, assume that the result is true for some $k \ge 2$. Let $m = k+1$ and $\displaystyle \beta \le 2^{\frac{1}{k+1}}$. Then certainly $\beta \le \displaystyle 2^{\frac{1}{k}}$, so by induction $T^i (a_2-a_1), \, T^i (\beta a_1-a_2) \in \Delta (0)$ for all $i \in \{ 0,1, \ldots, k-1 \}$. We only need to show that $T^k (a_2-a_1), \, T^k (\beta a_1-a_2) \in \Delta (0)$. First consider $T^k (a_2-a_1)= \beta^k (a_2-a_1)$. We have
\begin{eqnarray*}
T^k (a_2-a_1) \in \Delta (0) &\Leftrightarrow & \beta^k (a_2 -a_1) < \frac{a_1}{\beta}\\
&\Leftrightarrow & \frac{a_2}{\beta} < a_1 \, \left( \frac{1}{\beta} + \frac{1}{\beta^{k+2}} \right).
\end{eqnarray*}
Since $\displaystyle \beta \le 2^{\frac{1}{k+1}}$, then
$$ \frac{1}{\beta} + \frac{1}{\beta^{k+2}} \ge \frac{3}{2 \cdot 2^{\frac{1}{k+1}}} \ge \frac{3}{2\sqrt 2} >1.$$
Thus $T^k (a_2-a_1) \in \Delta (0)$. We now consider $T^k (\beta a_1-a_2) = \beta^k(\beta a_1 -a_2)$ and see that
\begin{eqnarray*}
T^k (\beta a_1-a_2) \in \Delta (0) &\Leftrightarrow & \beta^{k+1}a_1 - \beta^k a_2 < \frac{a_1}{\beta}\\
&\Leftrightarrow & a_1 \, \left( \beta - \frac{1}{\beta^{k+1}} \right) < a_2.
\end{eqnarray*}
Since $\displaystyle \beta \le 2^{\frac{1}{k+1}}$, then
$$\beta - \frac{1}{\beta^{k+1}} \le 2^{\frac{1}{k+1}} - \frac{1}{2} \le \sqrt 2 - \frac{1}{2} <1.$$
Thus $T^k (\beta a_1-a_2) \in \Delta (0)$ and this proves the lemma.
\end{proof}

\begin{remark}\label{r:smallbeta}
{\rm Suppose $2^{1/(m+1)} < \beta \le 2^{1/m}$ and $a_1 > \frac{a_2}{\beta}$. Lemma \ref{l:2to1overm} implies that $\kappa(i) = 2$ for $i \in \{ 1, \ldots, m\}$. By the same reasoning as in Remark \ref{r:wortel2}, $\kappa (n)$ would obtain the largest possible value if
$$\lambda (T^i (a_2-a_1) \cap \Delta(a_2))>0 \text{ and } \lambda (T^i(\beta a_1-a_2)  \cap \Delta(a_2))>0$$
for all $i= \ell m + (\ell-1)$, $\ell \ge 1$, and $T^i (a_2-a_1), T^i(\beta a_1-a_2) \in \Delta(0)$ for all other values of $i$. This would imply
$$
\bar{\kappa}(n) = \left\{
\begin{array}{ll}
\kappa(n), & \text{if $n = \ell m + (\ell -1)$ for some $\ell \ge 1$},\\
0, & \text{otherwise}.
\end{array}
\right.
$$
Thus, to get the maximal number of elements for $B_n$, we would have that if
$$(\ell -1)m +\ell \le n \le \ell m + \ell $$
for some $\ell$, then $\kappa (n) = 2^{\ell}$. So, in general we have that $\kappa (n) \le 2^{\lfloor \frac{n}{m} \rfloor +1 }$.
}
\end{remark}

For all $n \ge 1$, let $D_n$ be the union of all full fundamental intervals of rank $n$, that are not a subset of any full fundamental interval of lower rank. From the next lemma it follows that the full fundamental intervals generate the Borel $\sigma$-algebra on $[0, a_1)$.

\begin{lemma}\label{l:generate}
$$\lambda \left( \bigcup_{n \ge 1} D_n \right) = \sum_{n \ge 1} \lambda (D_n) = \lambda ([0, a_1)) = a_1.$$
\end{lemma}

\begin{proof}
Notice that all of the sets $D_n$ are disjoint. By (\ref{q:nonfull}) we have for each $n \ge 1$, that
$$ 0 \le \lambda \left( [0,a_1) \setminus \bigcup_{i=1}^n D_i \right) \le \kappa (n) \cdot \frac{a_1}{\beta^n},$$
so it is enough to prove that $\lim_{n \to \infty} \frac{\kappa(n)}{\beta^n} =0$. If $2 < \beta < 3$, then since $\kappa (n) \le 2^n$, we immediately have the result. For $1 < \beta \le G$, it follows from Remark \ref{r:wortel2} and Remark \ref{r:smallbeta}. Now, suppose $G < \beta \le 2$. Then by (\ref{q:Gbeta2}), we have that $\kappa (n) \le F(n+2)$, where $F(n+2)$ is the $(n+2)$-th element of the Fibonacci sequence. For the elements of this sequence, there is a closed formula which gives
\begin{equation}\label{q:fib}
F(n) = \frac{G^n - (1-G)^n}{\sqrt 5}.
\end{equation}
So
$$\frac{\kappa (n)}{\beta^n} \le \frac{1}{\sqrt 5} \left[ G^2 \left( \frac{G}{\beta} \right)^n - (1-G)^2 \left( \frac{1-G}{\beta} \right)^n \right].$$
Since $G < \beta \le 2$, also in this case $\lim_{n \to \infty} \frac{\kappa(n)}{\beta^n} =0$ and this proves the lemma.
\end{proof}

\begin{remark}\label{r:generate}
{\rm The fact that $\Delta(0)$ is a full fundamental interval of rank 1 allows us to construct full fundamental intervals of arbitrary small Lebesgue measure. This, together with the previous lemma, guarantees that we can write each subinterval of $[0,a_1)$ as a countable union of full fundamental intervals. Thus, the full fundamental intervals generate the Borel $\sigma$-algebra on $[0,a_1)$.
}
\end{remark}

Notice that for the cases illustrated by Figure \ref{f:attractorcases}(b) and \ref{f:attractorcases}(c), we can define the partitions $\Delta^{(n)}$, the sets $B_n$ and the numbers $\kappa (n)$ in a similar way. The only differences are that the support of the acim is given by the interval $[0, a_2-a_1)$ and that $\Delta(a_1)$ is the only full fundamental interval of rank 1. In that sense, $\Delta(a_1)$ plays the role of $\Delta(0)$ above. Since in these cases we have $2 < \beta < 3$ and since $\kappa (n) \le 2^n$ for all $n \ge 1$, we can prove a lemma similar to Lemma \ref{l:generate}, i.e. we can prove that the full fundamental intervals generate the Borel $\sigma$-algebra on the support of the acim.

\section{A natural extension of $T$}
For the version of the natural extension, we will define a space $R$, using the element of $B_n$. For $n \ge 1$, define the collections
$$R_n = \{ T^n \Delta(b_0 \ldots b_{n-1}) \times [0, \frac{a_1}{\beta^n} ) : \Delta(b_0 \ldots b_{n-1}) \in B_n \}.$$
So to each element of $B_n$, there corresponds an element of $R_n$ and thus the number of elements in $R_n$ is given by $\kappa (n)$. We enumerate the elements of $R_n$ and write $R_n = \{ R_{(n,i)} : 1 \le i \le \kappa (n)\}$. Thus, for each $\Delta(b_0 \ldots b_{n-1}) \in B_n$ there exists a unique $1 \le i \le \kappa (n)$ such that $T^n \Delta(b_0 \ldots b_{n-1}) \times [0, \frac{a_1}{\beta^n}) = R_{(n,i)}$. Let $R_0 = [0,a_1) \times [0,a_1)$ and let $R$ be the disjoint union of all these sets, i.e.
$$ R = R_0 \times \{0 \} \times \{0\} \cup \bigcup_{n=1}^{\infty} \bigcup_{i=1}^{\kappa(n)} R_{(n,i)} \times \{n \} \times \{i\}.$$
The $\sigma$-algebra $\mathcal R$ on $R$ is the disjoint union of the Borel $\sigma$-algebras on all the rectangles $R_{(n,i)}$ and $R_0$. Let $\lambda_R$ be the measure on $R$, given by the two dimensional Lebesgue measure on each of these rectangles. Define the probability measure $\nu$ on $R$ by setting $\nu (E) = \frac{1}{\lambda_R  (R)} \lambda_R (E)$ for all $E \in \mathcal R$. The next lemma says that this measure is well defined and finite.

\begin{lemma}\label{l:finitemeasure}
$$ \lambda_R(R) < \infty.$$
\end{lemma}

\begin{proof}
Set $\kappa (0)=1$. Then,
$$ \lambda_R (R) \le \sum_{n=0}^{\infty} \kappa (n) \frac{a_1^2}{\beta^n}.$$
Using the same arguments as in the proof of Lemma \ref{l:generate}, we can show that the sum on the right hand side converges for all $1< \beta <3$.
\end{proof}

\begin{figure}[h]
\centering
\subfigure[$\lambda(T^k (\beta a_1-a_2) \cap \Delta(a_1))=0$]{\includegraphics[width=.71\textwidth]{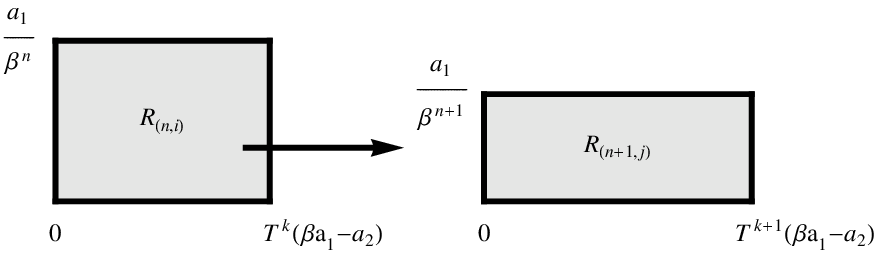}}
\quad
\subfigure[$\lambda(T^k (\beta a_1-a_2) \cap \Delta(a_1))>0$ and $\lambda(T^k (\beta a_1-a_2) \cap \Delta(a_2))=0$]{\includegraphics[width=.71\textwidth]{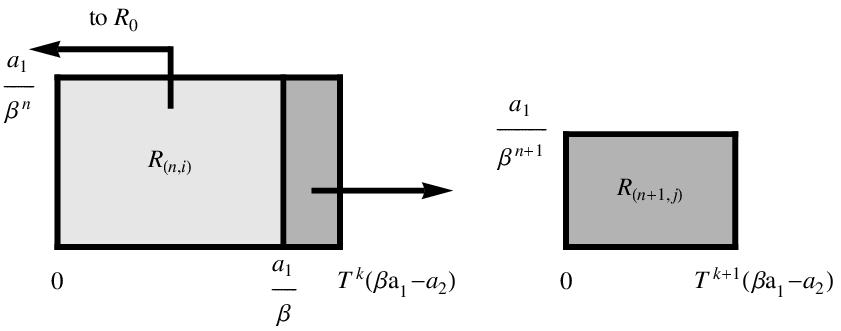}}
\quad
\subfigure[$\lambda(T^k (\beta a_1-a_2) \cap \Delta(a_2))>0$]{\includegraphics[width=.71\textwidth]{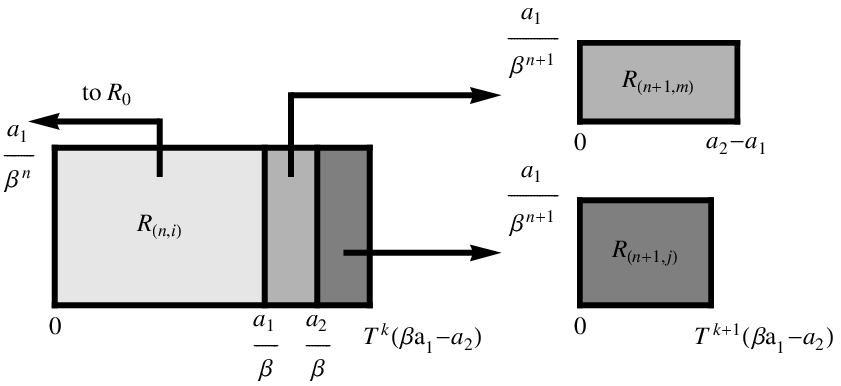}}
\caption{The arrows indicate the action of $\T$ in case $T^n \Delta(b_0 \ldots b_{n-1}) = T^k (\beta a_1 -a_2)$, for some $k < n$.}
\label{f:mathcalT}
\end{figure}

The transformation $\T: R \to R$ is defined piecewise on each rectangle. If $(x,y) \in R_0$, then
$$ \T (x,y,0,0) = \left\{
\begin{array}{ll}
(Tx,\frac{y}{\beta}, 0,0), & \text{if } x \in \Delta(0),\\
(Tx,\frac{y}{\beta},1,j), & \text{if } x \in \Delta(a_j), \; j\in \{ 1,2 \}.
\end{array}
\right.
$$
For each $n \ge 1$, $1 \le i \le \kappa(n)$, if
$$R_{(n,i)} = T^n \Delta(b_0 \ldots b_{n-1}) \times [0, \frac{a_1}{\beta^n} ),$$
then $\T$ maps this rectangle to the rectangles corresponding to the fundamental intervals contained in $\Delta(b_0 \ldots b_{n-1})$ in the following way. If $\Delta(b_0 \ldots b_{n-1}0)$ is full and $(x,y) \in R_{(n,i)}$ with $x \in \Delta (0)$, then
$$ \T (x,y,n,i) = (Tx, \frac{b_0}{\beta} + \frac{b_1}{\beta^2} + \ldots + \frac{b_{n-1}}{\beta^n} + \frac{y}{\beta},0,0).$$
If $\Delta(b_0 \ldots b_{n-1} b_n) \in B_{n+1}$ and $j$ is the index of the corresponding set in $R_{n+1}$, then for $(x,y) \in R_{(n,i)}$ with $x \in \Delta(b_n)$, we set
$$\T (x,y,n,i) = (Tx, \frac{y}{\beta}, n+1,j).$$
In Figure \ref{f:mathcalT} we show the different situations in case $T^n \Delta(b_0 \ldots b_{n-1}) = T^k (\beta a_1 -a_2)$ for some $k < n$. If $T^n \Delta (b_0 \ldots b_{n-1}) = T^k (a_2-a_1)$ for some $k<n$, the pictures look exactly the same with $T^k (a_2-a_1)$ in place of $T^k (\beta a_1 -a_2)$.

Notice that if a rectangle $R_{(n,i)}$ corresponds to a fundamental interval $\Delta (b_0 \ldots b_{n-1})$ such that $\Delta (b_0 \ldots b_{n-1}0)$ is non-full, then this is the only fundamental interval contained in $\Delta(b_0 \ldots b_{n-1})$. $\T$ then maps the rectangle $R_{(n,i)}$ bijectively to the rectangle $R_{(n+1, j)}$, corresponding to $\Delta(b_0 \ldots b_{n-1}0)$. Otherwise, $\Delta(b_0 \ldots b_{n-1}0)$ is a full fundamental interval contained in $\Delta(b_0 \ldots b_{n-1})$ and the rectangle $R_{(n,i)}$ is partly mapped surjectively onto some rectangle $R_{(n+1,j)}$ and partly into $R_0$. From Lemma \ref{l:generate} it now follows that $\T$ is bijective.\\
\vskip .3cm

Let $\pi_1:R \to [0, a_1)$ be the projection onto the first coordinate. Define the measure $\mu$ on $([0, a_1), \mathcal B([0, a_1)))$ by pulling back the measure $\nu$, i.e. for all measurable sets $E \in \mathcal B([0, a_1))$, let $\mu (E) = \nu (\pi_1^{-1}E)$. In order to show that $(R, \mathcal R, \nu, \T)$ is a version of the natural extension of the system $([0, a_1), \mathcal B([0, a_1)), \mu, T)$ with $\pi_1$ as the factor map, we will prove all of the following.
\begin{itemize}
\item[(i)] $\pi_1$ is a surjective, measurable and measure preserving map from $R$ to $[0,a_1)$.
\item[(ii)] For all $x \in R$, we have $(T \circ \pi_1)(x) = (\pi_1 \circ \T)(x)$.
\item[(iii)] $\T:R \to R$ is an invertible transformation.
\item[(iv)] $\mathcal R = \bigvee_{n=0}^{\infty} \T^n \pi_1^{-1} (\mathcal B([0,a_1)))$, where $\bigvee_{n=0}^{\infty} \T^n \pi_1^{-1} (\mathcal B([0,a_1)))$ is the smallest $\sigma$-algebra containing the $\sigma$-algebras $\T^n \pi_1^{-1} (\mathcal B ([0,a_1)))$ for all $n \ge 1$.
\end{itemize}
It is clear that $\pi_1$ is surjective and measurable and that $T \circ \pi_1 = \pi_1 \circ \T$. Since $\T$ expands by a factor $\beta$ in the first coordinate and contracts by a factor $\beta$ in the second coordinate, it is also clear that $\T$ is invariant with respect to the measure $\nu$. Then $\mu = \nu \circ \pi_1^{-1}$ defines a $T$-invariant probability measure on $([0,a_1), \mathcal B([0,a_1)))$, that is equivalent to the normalized Lebesgue measure on $[0,a_1)$ and $\pi_1$ is a measure preserving map. This shows (i) and (ii). The invertibility of $\T$ follows from Remark \ref{r:generate}, so that leaves only (iv). To prove (iv) we will have a look at the structure of the fundamental intervals and we will introduce some more notation.\\
\indent Let $\Delta (b_0 \ldots b_{n-1})$ be a fundamental interval. We can divide the block of digits $b_0 \ldots b_{n-1}$ into $M$ subblocks, $C_1, \ldots, C_M$, for some $M \ge 1$, where each subblock $C_i$, $1\le i \le M-1$, corresponds to a full fundamental interval. The last subblock, $C_M$, corresponds to a full fundamental interval exactly when $\Delta (b_0 \ldots b_{n-1})$ is full. We can make this precise, using the notion of return times to $R_0$. For points $(x,y) \in R_0$ define the {\it first return time to} $R_0$ by
$$r_1(x,y) = \inf \{ n \ge 1: \T^n (x,y,0,0) \in R_0 \times \{0\} \times \{0\} \}$$
and for $k \ge 1$, let the $k${\it -th return time to} $R_0$ be given recursively by
$$r_k (x,y) =\inf \{ n > r_{k-1} (x,y): \T^n (x,y,0,0) \in R_0 \times \{0\} \times \{0\} \}.$$
By the Poincar\'e Recurrence Theorem, we have $r_k (x,y) < \infty$ for almost all $(x,y) \in R_0$. Notice that this notion of return time depends only on $x$, i.e. for all $y, y' \in R_0$ and all $k \ge 1$, $r_k(x,y)=r_k (x,y')$. So we can write $r_k (x)$ instead of $r_k (x,y)$. If $\Delta (b_0 \ldots b_{n-1}) \in \Delta^{(n)}$, then for all $m \le n$, $\T^m$ maps the whole rectangle $\Delta(b_0 \ldots b_{n-1}) \times [0,a_1) \subseteq R_0$ to the same rectangle in $R$. So up to a certain $\ell \le n$, the $i$-th return time to $R_0$ is equal for all elements in $\Delta(b_0 \ldots b_{n-1})$. Now suppose that $\Delta(b_0 \ldots b_{n-1}) \in \Delta^{(n)}$ is a full fundamental interval, then there is an $M \ge 1$ and there are numbers $r_i$, $1 \le i \le M$ such that $r_i = r_i(x)$ for all $x \in \Delta (b_0 \ldots b_{n-1})$ and $r_M=n$. Put $r_0=0$. We can also obtain the numbers $r_i$ inductively as follows. Let
$$ r_1 = \inf\{ j >0 : T^{j+1} \Delta(b_0 \ldots b_j) = [0, a_1) \}$$
and if $r_1, \ldots, r_{k-1}$ are already known, let
$$r_k = \inf\{ j > r_{k-1} : T^{j+1} \Delta(b_{r_{k-1}} \ldots b_j) = [0, a_1) \}.$$
Take for $1 \le i \le M$, 
\begin{equation}\label{q:subblocks} 
C_i  = b_{r_{i-1}} \ldots b_{r_i-1}.
\end{equation}
Let $|C_i|$ denote the number of digits of the block $C_i$. The blocks have the following properties.
\begin{itemize}
\item[(p1)] For $1 \le i \le M$, $| C_i | = r_i-r_{i-1}$.
\item[(p2)] If $b_{r_i}=0$, then $r_{i+1} = r_i+1$. This means that if a subblock begins with the digit 0, then 0 is the only digit in this subblock. So, $C_{i+1}$ consists just of the digit 0.
\item[(p3)] For all $i \in \{ 1, \ldots, M \}$, $\Delta(C_i )$ is a full fundamental interval of rank $| C_i |$.
\end{itemize}

The next lemma is the last step in the proving that $(R, \mathcal R, \nu, \T)$ is the natural extension of the space $([0, a_1), \mathcal B([0, a_1)), \mu, T)$.

\begin{lemma}\label{l:bigvee}
Let $(R, \mathcal R, \nu, \T)$ and $([0, a_1), \mathcal B([0, a_1)), \mu, T)$ be the dynamical systems defined above. Then
$$\mathcal R = \bigvee_{n=0}^{\infty} \T^n (\pi_1^{-1} (\mathcal B([0, a_1)))).$$
\end{lemma}

\begin{proof}
It is clear that
$$\bigvee_{n=0}^{\infty} \T^n (\pi_1^{-1} (\mathcal B([0, a_1)))) \subseteq \mathcal R.$$
By Lemma \ref{l:generate} we know that the direct products of the full fundamental intervals contained in the rectangle $R_0$ generate the Borel $\sigma$-algebra on this rectangle. The same holds for all the rectangles $R_{(n,i)}$. First, let 
$$\Delta(d_0 \ldots d_{p-1}) \times \Delta(e_0 \ldots e_{q-1})$$
be a generating rectangle in $R_0$, where $\Delta(d_0 \ldots d_{p-1})$ and $\Delta(e_0 \ldots e_{q-1})$ are full fundamental intervals. For the set $\Delta(e_0 \ldots e_{q-1})$ construct the subblocks $ C_1 , \ldots, C_M$ as in (\ref{q:subblocks}). By property (p3) and Lemma \ref{l:fullfi}, $\Delta ( C_M C_{M-1} \ldots C_1  d_0 \ldots d_{p-1})$ is a full fundamental interval of rank $p+q$. Then
$$ \pi_1^{-1} (\Delta ( C_M  C_{M-1} \ldots C_1 d_0 \ldots d_{p-1})) \cap \left( R_0 \times \{0\} \times \{0\} \right) \quad$$
$$ \quad = \Delta ( C_M C_{M-1} \ldots C_1 d_0 \ldots d_{p-1}) \times [0,a_1) \times \{0 \} \times \{ 0 \}.$$
Since $\Delta(e_0 \ldots e_{q-1})$ is a full fundamental interval, it can be proven by induction that for all $i \in \{ 1, \ldots, q-1 \}$, $ T^i \Delta (e_0 \ldots e_{q-1}) = \Delta(e_i \ldots e_{q-1})$.
This, together with the definitions of the blocks $ C_i $ and the transformation $\T$ leads to
$$ \T^q ( \pi_1^{-1} ( C_M C_{M-1} \ldots C_1 d_0 \ldots d_{p-1})) \cap \left( R_0 \times \{0\} \times \{0\} \right)) \quad \quad$$
$$ \quad =  \Delta (d_0 \ldots d_{p-1}) \times \Delta ( C_1 C_2 \ldots C_M) \times \{ 0 \} \times \{ 0 \}.$$
So
$$ \Delta(d_0 \ldots d_{p-1}) \times \Delta(e_0 \ldots e_{q-1}) \times \{0 \} \times \{0 \} \subseteq \bigvee_{n=0}^{\infty} \T^n \pi_1^{-1}(\mathcal B([0,a_1))).$$
Now, for $n \ge 1$ and $i \in \{ 1,2, \ldots , \kappa(n) \}$, let $R_{(n,i)}$ be a rectangle in $R_n$ and suppose that it corresponds to the fundamental interval $\Delta(b_0 \ldots b_{n-1})\in B_n$. Hence, 
$$R_{(n,i)} = T^n \Delta(b_0 \ldots b_{n-1}) \times [0, \frac{a_1}{\beta^n}) \times \{ n \} \times \{ i \}.$$
Let $\Delta(d_0 \ldots d_{p-1}) \times \Delta(e_0 \ldots e_{q-1}) \times \{n \} \times \{i \}$ be a generating rectangle for the Borel $\sigma$-algebra on the rectangle $R_{(n,i)}$. So $\Delta(d_0 \ldots d_{p-1})$ and $\Delta(e_0 \ldots e_{q-1})$ are again full fundamental intervals. Notice that 
$$\Delta(e_0 \ldots e_{q-1}) \subseteq \Delta(\underbrace{0 \ldots 0}_{n \text{ times}}),$$
which means that $q \ge n$. Also, for all $i \in \{ 0, \ldots, n-1\}$, $e_i =0$ and thus $r_{i+1}=i+1$. So, if we divide $e_0 \dots e_{q-1}$ into subblocks $ C_i $ as before, we get that $ C_1 =  C_2 = \ldots = C_n  =0$, that $M \ge n$ and that $| C_{n+1}| + \ldots + | C_M|=q-n$. Consider the set 
$$ C = \Delta ( C_M C_{M-1} \ldots C_{n+1}  b_0 \ldots b_{n-1} d_0 \ldots d_{p-1}).
$$
We will show the following.\\
Claim: The set $C$ is a fundamental interval of rank $p+q$ and $T^q C=\Delta(d_0 \ldots d_{p-1})$.

\vskip .2cm

\noindent First notice that
$$ C= \Delta ( C_M C_{M-1} \ldots C_{n+1}) \cap T^{n-q}  \Delta(b_0 \ldots b_{n-1}) \cap T^{-q}\Delta(d_0 \ldots d_{p-1}).$$
So obviously,
$$
T^q C \subseteq T^q \Delta ( C_M C_{M-1} \ldots C_{n+1}) \cap T^n  \Delta(b_0 \ldots b_{n-1}) \cap \Delta(d_0 \ldots d_{p-1}).
$$
By Lemma \ref{l:fullfi}, $\Delta ( C_M C_{M-1} \ldots C_{n+1})$ is a full fundamental interval of rank $q-n$, so $ T^q \Delta ( C_M C_{M-1} \ldots C_{n+1}) = [0,a_1)$. Now, by the definition of $R_{(n,i)}$ we have that
\begin{equation}\label{q:deltaa}
\Delta(d_0 \ldots d_{p-1}) \subseteq T^n  \Delta(b_0 \ldots b_{n-1}),
\end{equation}
and thus $T^q C \subseteq   \Delta(d_0 \ldots d_{p-1})$.\\
For the other inclusion, let $z \in \Delta (d_0 \ldots d_{p-1})$. By (\ref{q:deltaa}),
there is an element $y$ in $\Delta(b_0 \ldots b_{n-1})$, such that $T^n y =z$. And since $ T^{q-n} \Delta ( C_M C_{M-1} \ldots C_{n+1}) = [0,a_1)$, there is an $x \in \Delta ( C_M C_{M-1} \ldots C_{n+1})$ with $T^{q-n}x = y$, so $T^q x = z$. This means that
$$ z \in T^q \Delta ( C_M C_{M-1} \ldots C_{n+1}) \cap T^n \Delta(b_0 \ldots b_{n-1}) \cap \Delta(d_0 \ldots d_{p-1}).$$ 
So $ T^q C= \Delta(d_0 \ldots d_{p-1})$ and this proves the claim.

\vskip .2cm

\noindent Consider the set $D= \pi_1^{-1} (C) \cap \left( R_0 \times \{ 0 \} \times \{0 \} \right)$. Then as before, we have
$$\T^{q-n} D = \Delta (b_0 \ldots b_{n-1}d_0 \ldots d_{p-1}) \times \Delta( C_{n+1} C_{n+2}\ldots C_M) \times \{0 \} \times \{0 \}.$$
And after $n$ more steps,
\begin{eqnarray*}
\T^q D &=& \Delta(d_0 \ldots d_{p-1}) \times \Delta (\underbrace{00\ldots 0}_{n \text{ times}} C_{n+1} \ldots C_M) \times \{ n \} \times \{ i \}\\
&=& \Delta(d_0 \ldots d_{p-1}) \times \Delta(e_0 \ldots e_{q-1}) \times \{n\} \times \{i \}.
\end{eqnarray*}
So, 
$$\Delta(d_0 \ldots d_{p-1}) \times \Delta(e_0 \ldots e_{q-1}) \times \{n \} \times \{ i \} \in \bigvee_{n=0}^{\infty} \T^n \pi_1^{-1}(\mathcal B([0,a_1)))$$
and thus we see that
$$ \mathcal R = \bigvee_{n=0}^{\infty} \T^n \pi_1^{-1}(\mathcal B([0,a_1))). \qedhere$$
\end{proof}

This gives the following theorem.
\begin{theorem}
The dynamical system $(R, \mathcal R, \nu, \T)$ is a version of the natural extension of the dynamical system $([0,a_1), \mathcal B([0,a_1)), \mu, T)$. Here, the measure $\mu$, given by $\mu(E) = \nu (\pi_1^{-1}(E))$ for all measurable sets $E$, is the acim of $T$ and the density of $\mu$ is equal to the density from equation (\ref{q:wilkinson}).
\end{theorem}

\begin{proof}
The fact that $(R, \mathcal R, \nu, \T)$ is a version of the natural extension of the system $([0,a_1), \mathcal B([0,a_1)), \mu, T)$ follows from Remark \ref{r:generate}, the properties of the map $\pi_1$ and Lemma \ref{l:bigvee}. Now for each measurable set $E \in \mathcal B ([0,a_1))$, we have
\begin{eqnarray*}
\mu (E) &=& \nu (\pi_1^{-1}(E))\\
&=& \frac{1}{\lambda_R(R)} \left[ a_1 \lambda(E) + \sum_{n=1}^{\infty} \sum_{\Delta(b_0 \ldots b_{n-1}) \in B_n} \frac{a_1}{\beta^n} \lambda (E \cap T^n \Delta(b_0 \ldots b_{n-1}))\right]\\
&=& \frac{a_1}{\lambda_R(R)} \int_E \left[ 1+ \sum_{n=1}^{\infty} \sum_{\Delta(b_0 \ldots b_{n-1}) \in B_n} \frac{a_1}{\beta^n} \; 1_{T^n \Delta(b_0 \ldots b_{n-1})} \right] d\lambda.
\end{eqnarray*}
Here $T^n \Delta(b_0 \ldots b_{n-1})$ has the form
$$[0, T^i (a_2-a_1)) \text{ or } [0, T^i (\beta a_1 - a_2))$$
for some $0 \le i <n$. So, the density of the measure $\mu$ equals the density from equation (\ref{q:wilkinson}).
\end{proof}

\begin{remark}
{\rm (i) The above definitions of the space $R$ and the transformation $\T$ can be adapted quite easily for the cases illustrated by Figure \ref{f:attractorcases}(b) and \ref{f:attractorcases}(c). We let $\Delta(a_1)$ take the role of $\Delta(0)$ and consider the orbits of the points $a_1$ and $\beta (a_2-a_1)-a_2$. In general, the sets $B_n$ will contain more elements, but since $2 < \beta < 3$, it is immediate that $\lambda_R(R) < \infty$. This shows that we can construct a version of the natural extension of $T$, also for these two cases.\\
(ii) Let $R_0'$ be the set obtained from $R_0$ by removing the set of measure zero of elements which do not return to $R_0$, i.e. we remove those $(x,y)$ for which $r_1(x,y)=\infty$. Let $\mathcal W :R_0' \to R_0'$ be the transformation induced by $\T$, i.e. for all $(x,y) \in R_0'$, let
$$\mathcal W (x,y) = \T^{r_1(x,y)}(x,y,0,0).$$
Then the dynamical system $(R_0', \mathcal B(R_0'), \lambda \times \lambda, \mathcal W)$, where $\mathcal B(R_0')$ is the Borel $\sigma$-algebra on $R_0'$, is isomorphic to the natural extension of a GLS-transformation as defined in \cite{DKS}. This implies that the system $(R_0', \mathcal B(R_0'), \lambda \times \lambda, \mathcal W)$ is Bernoulli.
}
\end{remark}

Using this invariant measure, we will show that $T$ is an exact transformation. Since the full fundamental intervals generate the Borel $\sigma$-algebra on the support of the acim, by a result of Rohlin (\cite{Roh1}), it is enough to show that there exists a universal constant $\gamma >0$, such that for any full fundamental interval $\Delta(b_0 \ldots b_{n-1})$ and any measurable subset $E \subseteq \Delta(b_0 \ldots b_{n-1})$, we have
$$ \mu (T^n E) \le \gamma \cdot \frac{\mu(E)}{\mu(\Delta(b_0 \ldots b_{n-1}))}.$$
To this end, define two constants, $c_1,c_2 > 0$, by
$$c_1 = \frac{a_1}{\lambda_R (R)}  \quad \text{and} \quad c_2 = 1+\sum_{n=1}^{\infty} \sum_{\Delta(b_0 \ldots b_{n-1})\in B_n} \frac{1}{\beta^n}.$$
Then for all measurable sets $E$, we have
$$c_1 \lambda (E) \le \mu(E) \le c_1 c_2 \lambda(E).$$
Now, let $\Delta (b_0 \ldots b_{n-1})$ be a full fundamental interval of rank $n$. Then by (\ref{q:genfull}),
$$\lambda (\Delta (b_0 \ldots b_{n-1})) = \frac{a_1}{\beta^n}.$$
Let $E \subseteq \Delta (b_0 \ldots b_{n-1})$ be a measurable set. Then
$$\lambda (T^n E) = \beta ^n \lambda (E) = \frac{a_1}{\lambda (\Delta (b_0 \ldots b_{n-1}))} \lambda (E).$$
Now,
\begin{eqnarray*}
\mu (T^n E) &\le& c_1 c_2 \lambda (T^n E) = c_1c_2 a_1 \frac{\lambda (E)}{\lambda (\Delta (b_0 \ldots b_{n-1}))}\\
& \le & c_1 c_2 a_1 \frac{\mu(E) c_1 c_2}{c_1 \mu(\Delta (b_0 \ldots b_{n-1}))} =  c_1 c_2^2 a_1 \frac{\mu(E)}{\mu(\Delta (b_0 \ldots b_{n-1}))}.
\end{eqnarray*}
If we take $\gamma = c_1 c_2^2 a_1$, then $\gamma >0$ and
$$ \mu (T^n E) \le \gamma \cdot \frac{\mu(E)}{\mu(\Delta(b_0 \ldots b_{n-1}))}.$$
Thus, $T$ is exact and hence mixing of all orders. Furthermore, the natural extension $\T$ is a $K$-automorphism. By a result of Rychlik (see \cite{Ryc1}), it follows immediately that $T$ is weakly Bernoulli.

\section{A specific example}
\begin{figure}[h]
\centering
{\includegraphics[height=5cm]{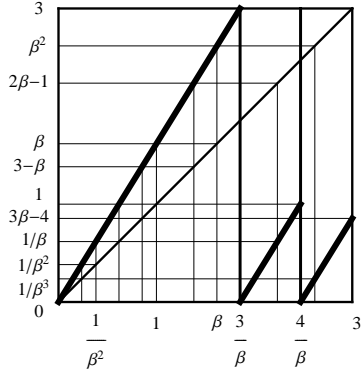}}
\caption{The orbits of $1$ and $3\beta-4$ under the greedy $\beta$-transformation with $\beta = \frac{1+\sqrt 5}{2}$ and $A= \{ 0,3,4 \}$.}
\label{f:goud034}
\end{figure}

In this section we will consider one specific example of a $\beta$-transformation with three deleted digits. Let $G$ be the golden mean as before and take $\beta = G$. Consider the allowable digit set $A= \{ 0,3,4 \}$. The acim of the greedy $\beta$-transformation with deleted digits with this $\beta$ and $A$, has the interval $[0,3)$ as its support. The partition $\Delta = \{ \Delta(0), \Delta(3), \Delta(4) \}$ is given by
$$\Delta(0) = \left[ 0, \frac{3}{\beta} \right), \quad \Delta(3)= \left[ \frac{3}{\beta}, \frac{4}{\beta} \right), \quad \Delta(4)=\left[ \frac{4}{\beta}, 3 \right) $$
and the transformation is then $Tx = \beta x -j$, if $x \in \Delta(j)$. The points $x=1$ and $x=3\beta -4$ are of special interest and their orbits under $T$ are as follows.
$$
\begin{array}{llll}
1, & T1= \beta, & T^2 1 =  \beta^2, & T^3 1 = \displaystyle \frac{1}{\beta^3},\\
& T^4 1=  \frac{1}{\beta^2}, & T^5 1 =  \frac{1}{\beta}, & T^6 1 =1,\\
3\beta -4, & T (3\beta -4) = 3-\beta, & T^2 (3\beta-4) = 2\beta-1 & T^3 (3\beta-4) = \displaystyle \frac{1}{\beta}.
\end{array}
$$
Thus their greedy $\beta$-expansions are given by
\begin{eqnarray*}
1 &=& \displaystyle \sum_{n=1}^{\infty} \frac{d_n^{(1)}}{\beta^n} =_{\beta} \overline{00400},\\
3\beta -4 &=& \displaystyle \sum_{n=1}^{\infty} \frac{d_n^{(2)}}{\beta^n} =_{\beta} 003\overline{00004},\\
\end{eqnarray*}
where the bar on the right hand side of these equations indicates a repeating block in the expansions. Notice that $3\beta -4$ would be the image of 3 under $T$, if $T$ were not restricted to the half open interval $[0,3)$. Figure \ref{f:goud034} shows the orbits of both points under $T$.

\begin{figure}[h]
\centering
\includegraphics[width=\textwidth]{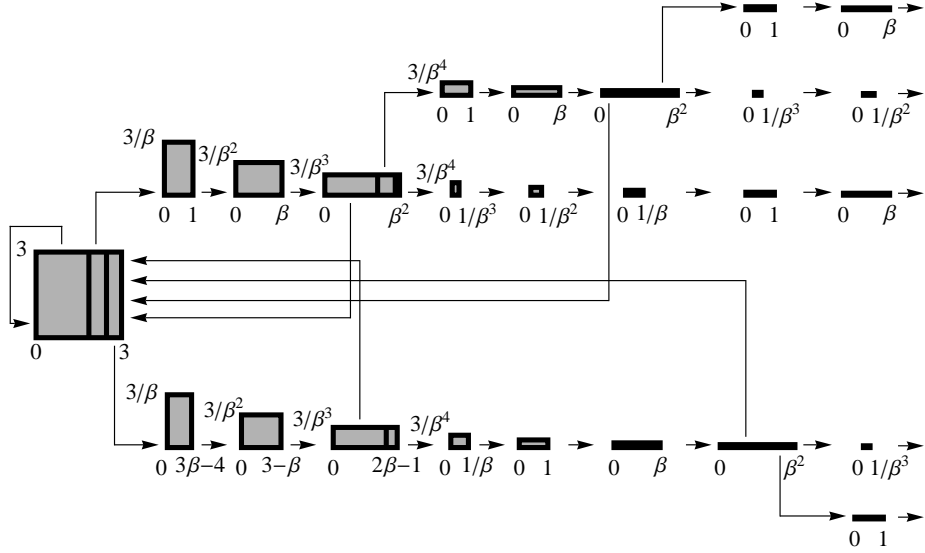}
\caption{The space $R$ consists of all these rectangles. The arrows indicate where the rectangles are mapped under $\T$.}
\label{f:tree}
\end{figure}

For each non-full fundamental interval, $\Delta(b_0 \ldots b_{n-1})$, the set $T^n \Delta(b_0 \ldots b_{n-1})$ is one of the following,
$$
\begin{array}{lllll}
[0, 3\beta-4), & [0, 3-\beta), & [0,2\beta-1), & [0,1), & [0,\beta),\\
&[0, \beta^2), & [0, 1/ \beta^3), & [0, 1/ \beta^2 ), & [0, 1/ \beta). 
\end{array}
$$
This means that we could give all the elements of $R_n$ explicitly. We will not do this, since Figure \ref{f:tree} speaks for itself. The space $R$ contains all the rectangles shown in the figure.\\
\indent To determine the invariant measure of $T$, we need to project the Lebesgue measure on the rectangles of $R$ onto the first coordinate. Therefore, we need to add up the heights of all the rectangles which have the same interval in the first coordinate. To do this, we only need to determine the total height of the rectangles of the form $[0,1) \times [0, \frac{3}{\beta^n})$, since the total height of all the other rectangles in $R$ can be deduced from this. These heights can be found by using the Fibonacci numbers $F(n)$, as defined before. Then $\kappa(1)=2$ and observe that for $n \ge 2$, $\kappa (n) = F(\lfloor \frac{n-1}{3} \rfloor +2) + F(\lfloor \frac{n-2}{3} \rfloor+1)$. For $n=3k+1$, $k \ge 0$, the number of rectangles in $R$ of the form $[0,1) \times [0, \frac{3}{\beta^n})$ is equal to $F(k+1)$ and for $n = 3k+2$, $k \ge 1$, this number is equal to $F(k)$. Formula (\ref{q:fib}) gives that the total Lebesgue measure of all these rectangle is equal to
\begin{eqnarray*}
\sum_{k=0}^{\infty} \frac{3F(k+1)}{\beta^{3k+1}} + \sum_{k=1}^{\infty} \frac{3F(k)}{\beta^{3k+2}} &=& \frac{3}{\beta}\left[1+ \sum_{k=1}^{\infty} \frac{\beta F(k+1) +F(k)}{\beta^{3k+1}} \right]\\
&=&  \frac{3}{\beta}\left[1+\frac{1}{\sqrt 5}\left(\beta \sum_{k=1}^{\infty} \frac{1}{\beta^{2k}} + \frac{1}{\beta}\sum_{k=1}^{\infty} \frac{1}{\beta^{2k}} \right)\right]\\
&=&  \frac{3}{\beta}\left[1+\frac{1}{\sqrt 5}\left( (\beta + \frac{1}{\beta}) (\frac{\beta^2}{\beta^2-1} -1 ) \right)\right]\\
&=&  \frac{3}{\beta}\left[1+\frac{1}{\sqrt 5}(3-\beta)\right]=3.
\end{eqnarray*}

The total height of the rectangles of the form $[0, \beta) \times [0, \frac{3}{\beta^n})$ is now equal to $\frac{3}{\beta}$, that of the rectangles of the form $[0, \beta^2) \times [0, \frac{3}{\beta^n})$ is $\frac{3}{\beta^2}$, etcetera. The total height of the rectangles $[0, \frac{1}{\beta} ) \times [0, \frac{3}{\beta^n})$ is given by $\frac{3}{\beta^4}+\frac{3}{\beta^5} = \frac{3}{\beta^3}$. Then the density function of the invariant probability measure of $T$, $h:[0,3) \to [0,3)$, equivalent to the normalized Lebesgue measure on $[0,3)$, is given by
\begin{eqnarray*}
h(x) &=&\frac{1}{27-4\beta} \left[ \frac{1}{\beta} \; 1_{[0, 3\beta-4)}(x) + \frac{1}{\beta^2} \; 1_{[0,3-\beta)}(x) + \frac{1}{\beta^3} \; 1_{[0, 2\beta-1)}(x) \right. \\
&& + 1_{[0,1)}(x) + \frac{1}{\beta} \; 1_{[0, \beta)}(x) + \frac{1}{\beta^2} \; 1_{[0, \beta^2)}(x) + \frac{1}{\beta^3} \; 1_{[0, 1/ \beta^3)}(x)\\
&& \left. + \frac{1}{\beta^4} \; 1_{[0, 1 / \beta^2)}(x) + \frac{1}{\beta^3} 1_{[0, 1/ \beta)}(x) + 1_{[0,3)}(x) \right] .
\end{eqnarray*}

\section{Concluding remarks}
We have seen that it is possible to find an explicit expression for the density function of the acim of the greedy $\beta$-transformation with three deleted digits. This density function is equal to the density found by Wilkinson in \cite{Wil1}. We have constructed a version of the natural extension of the dynamical system of the greedy $\beta$-transformation with deleted digits to find this density in case the real number $\beta>1 $ and the digit set $A$ satisfy condition (\ref{q:conditions}). A very similar construction can be used to give a version of the natural extension for the transformation in case $A= \{ 0, a_1, a_2 \}$ does not satisfy this condition. The transformation $T$ is exact and weakly Bermoulli.\\
\indent Of course, it still remains to find an explicit expression for the acim of the greedy $\beta$-transformation with deleted digits, of which the digit set contains more than three digits. If the number of digits, $m$, satisfies $m < \beta \le m+1$, the density is given by the density from (\ref{q:wilkinson}). For this case, one could try to generalize the definition of the version of the natural extension that we gave in this paper. The space $R$ would then become an $m$-ary tree and the finiteness of the measure of $R$ is guaranteed, since $\beta >m$.\\
Recently, G\'ora found an invariant density for piecewise linear and increasing maps of constant slope, for certain values of $\beta$. See \cite{Gor1}.

\def\cprime{$'$}

\end{document}